\documentclass[12pt]{amsart}

\usepackage{amssymb,latexsym, fullpage, eucal}
\usepackage{enumerate}

\theoremstyle{definition}

\numberwithin{equation}{section}

\begin{document}


%
%

\title[A note on a conjecture of Gonek]{A note on a conjecture of Gonek}

\author[Micah B. Milinovich]{Micah B. Milinovich}
\address{Department of Mathematics, University of Mississippi, University, MS  38677 USA}
\email{mbmilino@olemiss.edu}
\author[Nathan Ng]{Nathan Ng}
\address{University of Lethbridge, Department of Mathematics and Computer Science, 4401 University Drive, Lethbridge, AB Canada T1K 3M4 }
\email{nathan.ng@uleth.ca}

\thanks{MBM is supported in part by a University of Mississippi College of Liberal Arts summer research grant.  NN is supported in part by an NSERC Discovery grant.}




\maketitle


\begin{abstract}
We derive a lower bound for a second moment of the reciprocal of the derivative of the Riemann zeta-function
over the zeros of $\zeta(s)$  that is half the size of the conjectured value.  Our result is conditional upon the assumption of the Riemann Hypothesis and the conjecture that the zeros of the zeta-function are simple.
\end{abstract}


\maketitle

\section{Introduction}

Let $\zeta(s)$ denote the Riemann zeta-function.  Using a heuristic method similar to Montgomery's study \cite{M} of the pair-correlation of the imaginary parts of the non-trivial zeros of $\zeta(s)$,  Gonek has made the following conjecture  \cite{G1, G3}. \\

\noindent{\bf Conjecture.} {\it Assume the Riemann Hypothesis and that the zeros of $\zeta(s)$ are simple.  Then, as $T\rightarrow \infty$,}
\begin{equation}\label{conj}
\sum_{0<\gamma\leq T} \frac{1}{{|\zeta'(\rho)|}^{2}} \sim \frac{3}{\pi^{3}} T
\end{equation}
{\it where the sum runs over the non-trivial zeros $\rho=\tfrac{1}{2}\!+\!i\gamma$ of $\zeta(s)$.}\\

\noindent The assumption on the simplicity of the zeros of the zeta-function in the above conjecture is so that the sum over zeros on the right-hand side of (\ref{conj}) is well defined.  While the details of Gonek's method have never been published, he announced his conjecture in \cite{G2}. More recently, a different heuristic method of Hughes, Keating, and O'Connell \cite{HKO} based upon modeling the Riemann zeta-function and its derivative using the characteristic polynomials of random matrices has led to the same conjecture.  Through the work of Ingham \cite{In}, Titchmarsh (Chapter 14 of \cite{T}), Odlyzko and te Riele \cite{OR}, Gonek (unpublished), and Ng \cite{N}, it is known that the behavior of this and related sums are intimately connected to the distribution of the summatory function
\begin{equation*}
M(x) = \sum_{n\leq x} \mu(n)
\end{equation*}
where $\mu(\cdot)$, the M\"{o}bius function, is defined by $\mu(1)=1$, $\mu(n)=(-1)^k$ if $n$ is divisible by $k$ distinct primes, and $\mu(n)=0$ if $n>1$ is not square-free.  See also \cite{H} and \cite{St} for connections between similar sums and other arithmetic problems.

In support of his conjecture, Gonek \cite{G2} has shown, assuming the Riemann Hypothesis and the simplicity of the zeros of $\zeta(s)$, that
\begin{equation}\label{neg1}
  \sum_{0<\gamma\leq T} \frac{1}{{|\zeta'(\rho)|}^{2}} \geq C T
 \end{equation}
for some constant $C>0$ and $T$ sufficiently large. In this note, we show that the inequality in (\ref{neg1}) holds for any constant $C < \frac{3}{2\pi^3}$. \\

\noindent{\bf Theorem. }{\it  Assume the Riemann Hypothesis and that the zeros of $\zeta(s)$ are simple.  Then,
for any fixed $\varepsilon >0$,}
\begin{equation}\label{neg2}
\sum_{0<\gamma\leq T} \frac{1}{{|\zeta'(\rho)|}^{2}} \geq \left( \frac{3}{2\pi^{3}} - \varepsilon\right) T
\end{equation}
{\it for $T$ sufficiently large.}\\

While our result differs from the conjectural lower bound by a factor of $2$, any improvements in the strength of this lower bound have, thus far, eluded us.  
It would be interesting to investigate whether for $k >0$ there  is a constant $C_k>0$ such that
\begin{equation}\label{neg}
\sum_{0<\gamma\leq T} \frac{1}{{|\zeta'(\rho)|}^{2k}} \geq C_k T (\log T)^{(k-1)^2}
\end{equation}
for $T$  sufficiently large.  However, a lower bound of this form is probably not of the correct order of magnitude for all $k$.  This is because it is expected that for each $\varepsilon>0$
there are infinitely many zeros $\rho=\tfrac{1}{2}+i\gamma$ of $\zeta(s)$ satisfying $|\zeta'(\rho)|^{-1} \gg |\gamma|^{1/3-\varepsilon}$.  
If such a sequence were to exist, it would then follow that
$$ \sum_{0<\gamma\leq T} \frac{1}{{|\zeta'(\rho)|}^{2k}} = \Omega\left( T^{2k/3-\varepsilon} \right)$$
and the lower bound in (\ref{neg}) would be significantly weaker than this $\Omega$-result when $k > \frac{3}{2}$.


\section{Proof of Theorem}

The method we use to prove our theorem is based on a recent idea of Rudnick and Soundararajan \cite{RS1}. Let
\begin{equation}  \label{eq:Tth}
   \xi\!=\!T^\vartheta \end{equation}
where $0 < \vartheta <1$ is fixed and define the Dirichlet polynomial
\begin{equation*}
\mathcal{M}_{\xi}(s)\!=\!\sum_{n\leq \xi} \mu(n) n^{-s}
\end{equation*}
where $\mu$ is the M\"{o}bius function.   Assuming the Riemann Hypothesis, for any non-trivial zero $\rho=\tfrac{1}{2}+i\gamma$ of $\zeta(s)$, we see that $\overline{\mathcal{M}_{\xi}(\rho)} = \mathcal{M}_{\xi}(1\!-\!\rho)$.  From
 this observation and Cauchy's inequality it follows that
\begin{equation}\label{eq:cauchy}
 \sum_{0<\gamma\leq T} \frac{1}{{|\zeta'(\rho)|}^{2}} \geq \frac{ \big|M_{1}\big|^{2}}{M_{2}}
 \end{equation}
where
 $$ M_{1} =  \sum_{0<\gamma\leq T} \frac{1}{\zeta'(\rho)} \mathcal{M}_{\xi}(1\!-\!\rho) \quad \text{ and } \quad M_{2} =\sum_{0<\gamma\leq T} \big| \mathcal{M}_{\xi}(\rho)\big|^{2}.$$
Our Theorem is a consequence of the following proposition. \\

\noindent{\bf Proposition.} {\it Assume the Riemann Hypothesis and let $0<\vartheta<1$ be fixed. Then}
\begin{equation}
\label{eq:m2}
M_2= \frac{3}{\pi^{3}} \left( \vartheta+\vartheta^2\right) T \log^2 T+ O(T \log T).
\end{equation}
{\it If we further assume that the zeros of  $\zeta(s)$ are all simple, then there exists a sequence $\mathcal{T}:=\{ \tau_n \}_{n =3}^{\infty}$ such that $n < \tau_n \le n+1$ and for $T \in \mathcal{T}$ we have}
\begin{equation}
\label{eq:m1}
M_1= \frac{3 \vartheta}{\pi^{3}} T \log T+ O(T).
\end{equation}

\bigskip

We now deduce our theorem from the above proposition.

\begin{proof}[Proof of the Theorem]

Let $T \geq 4$ and choose $\tau_n$ to satisfy $T-1 \le \tau_n < T$.
Combining \eqref{eq:cauchy}, \eqref{eq:m1}, and \eqref{eq:m2} we see that 
\begin{equation}
\begin{split}
\label{ineqtheta}
\sum_{0<\gamma\leq T} \frac{1}{{|\zeta'(\rho)|}^{2}} \geq 
 \sum_{0 <\gamma\leq \tau_n} \frac{1}{{|\zeta'(\rho)|}^{2}} &   \ge
\frac{\vartheta^2}{(\vartheta+\vartheta^2)} \frac{3}{\pi^{3}} \ \! \tau_n + o(\tau_n) \\
& \ge \frac{1}{(1+\vartheta^{-1})} \frac{3}{\pi^{3}} \ \! T + o(T)
\end{split}
\end{equation}
under the assumption of the Riemann Hypothesis and the simplicity of the zeros of $\zeta(s)$.  From (\ref{ineqtheta}),  our theorem follows by letting $\vartheta\rightarrow 1^-$.
\end{proof}

We could have just as easily estimated the sums $M_1$ and $M_2$ using a Dirichlet polynomial $\sum_{n \le \xi} a_n n^{-s}$ for a large class of coefficients $a_n$ in place of $\mathcal{M}_\xi(s)$. In the special case where $$a_n = \mu(n) P\big(\tfrac{\log \xi/n}{\log \xi}\big)$$ for polynomials $P$, we can show that the choice $P=1$ is optimal in the sense that it leads to largest lower bound in (\ref{neg2}).

We prove the above proposition in the next two sections; the sum $M_1$ is estimated in section 3 and  the sum $M_2$ is estimated in section 4. 
The evaluation of sums like $M_1$ dates back to Ingham's \cite{In} important work on $M(x)$ in which he considered sums of the form 
\[
  \sum_{0 < \gamma < T} (T-\gamma)^k \zeta'(\rho)^{-1}  
\]
for $k \in \mathbb{R}.$ The sum $M_2$  is of the form
\begin{equation}
  \label{eq:dism}
\sum_{0 < \gamma < T}|A(\rho)|^2 \quad \text{ where } \quad  A(s)= \sum_{n \le \xi} a_n n^{-s}
\end{equation}
is a Dirichlet polynomial with $\xi\leq T$. Such sums have played an important role in various applications. 
For instance,
results concerning the distribution of consecutive zeros of $\zeta(s)$ and discrete mean values of the zeta-function and its derivatives are proven  in  \cite{BMN,CGG,Cfivea,G5,Mi,N2,So}.
In each of these articles, the evaluation of the discrete mean \eqref{eq:dism} either makes use of the Guinand-Weil explicit formula or of Gonek's uniform version \cite{G5} of  Landau's formula
\begin{equation}
 \label{eq:langon}
\sum_{\substack{0 < \gamma < T \\ \zeta(\beta+i \gamma)=0}} 
 x^{\beta+i \gamma} = -\frac{T}{2\pi}\Lambda(x) +E(x,T)
\end{equation}
for $x,T > 1$ where $E(x,T)$ is an explicit error function uniform in $x$ and $T$.  A novel aspect of our approach is that it does not require the use of the Guinand-Weil explicit formula or of the Landau-Gonek explicit formula \eqref{eq:langon}. Instead we evaluate $M_2$ using the residue theorem and a version of Montgomery and Vaughan's mean value theorem for Dirichlet polynomials \cite{MV}. Our approach is simpler and it is likely that it can be extended to evaluate the discrete mean \eqref{eq:dism} for a large class of coefficients $a_n$ with $\xi \le T$. \\


\section{The estimation of $M_{1}$}

To estimate $M_{1},$  we require the following version of  Montgomery and Vaughan's mean value theorem for Dirichlet polynomials. \\

\noindent{\bf Lemma.} {\it Let $\{a_n\}$ and $\{b_n\}$ be two sequences of complex numbers. For any real number $T>0$, we have}
\begin{equation}
\begin{split}
   \label{eq:mv}
     \int_{0}^{T}  
   \Bigg(
   \sum_{n=1}^{\infty} a_n n^{-it}  \Bigg)&
    \Bigg( \sum_{n=1}^{\infty} b_n n^{it}  \Bigg) 
    dt  
    \\&= T\sum_{n=1}^{\infty} a_n b_n 
    + O \Bigg( \Big( 
   \sum_{n=1}^{\infty} n |a_n|^2
   \Big)^{\frac{1}{2}} 
    \Big( 
   \sum_{n=1}^{\infty} n |b_n|^2
   \Big)^{\frac{1}{2}} 
      \Bigg).
\end{split}
\end{equation}

\begin{proof}
This is Lemma 1 of Tsang \cite{Ts}. The special case where $b_n=\overline{a_n}$, is originally due to Montgomery and Vaughan \cite{MV}. It turns out, as shown by Tsang, that this special case is equivalent to the more general case stated in the lemma.
\end{proof}

Let $T \geq 4$ and set $c = 1+ (\log T)^{-1}$. 
It is well known (see Theorem 14.16 of Titchmarsh \cite{T}) that assuming the Riemann Hypothesis there exists a sequence $\mathcal{T} = \{\tau_{n}\}_{n=3}^{\infty}, \ n< \tau_{n} \leq n+1$, and a fixed constant $A>0$ such that
\begin{equation}\label{neg3}
\big|\zeta(\sigma\!+\! i\tau_{n})\big|^{-1} \ll \exp\Big(\frac{ A \log\tau_{n}}{\log\log\tau_{n}}\Big)
\end{equation}
uniformly for $\tfrac{1}{2}\leq \sigma \leq 2.$  We now prove the estimate \eqref{eq:m1} assuming that $T \in \mathcal{T}$.
Recall that $|\gamma|>1$ for every non-trivial zero $\rho=\tfrac{1}{2}+i\gamma$ of $\zeta(s)$.  Thus, assuming that all the zeros of $\zeta(s)$ are simple, the residue theorem implies that
\begin{equation*}
\begin{split}
 M_{1} &= \frac{1}{2\pi i} \left(\int_{c+i}^{c+iT} + \int_{c+iT}^{1-c+iT}+\int_{1-c+iT}^{1-c+1}+\int_{1-c+i}^{c+i}\right) \ \frac{1}{\zeta(s)} \mathcal{M}_{\xi}(1\!-\!s) \ ds 
 \\
 & = I_{1}+I_{2}+I_{3}+I_{4},
 \end{split}
 \end{equation*}
say.  Here we are using the fact that the residue of the function $1/\zeta(s)$ at $s=\rho$ equals $1/\zeta'(\rho)$ if $\rho$ is a simple zero of $\zeta(s)$.

The main contribution to $M_{1}$ comes from the integral $I_{1}$; the remainder of the integrals contribute an error term.   Observe that 
\[
  I_1 = \frac{1}{2 \pi} \int_{1}^{T} \sum_{m=1}^{\infty} \frac{\mu(m)}{m^{c+it}}   \sum_{n \le \xi } \frac{\mu(n)}{n^{1-c-it}} dt. 
\]
By \eqref{eq:mv} with $a_m = \mu(m)m^{-c}$ and $b_n = \mu(n)n^{-1+c}$ it follows that 
\begin{align*}
   I_1 & =  \frac{(T-1)}{2 \pi} \sum_{n \le \xi} \frac{\mu(n)^2}{n} + O \Bigg(
  \Big(  \sum_{n=1}^{\infty}  \frac{\mu(n)^2}{n^{2c-1}} \Big)^{\frac{1}{2}}
    \Big(  \sum_{n \le \xi} \mu(n)^2 n^{2c-1} \Big)^{\frac{1}{2}}
   \Bigg).
\end{align*}
Since 
\begin{equation}
 \label{eq:neg4}
 \sum_{n\leq \xi} \frac{\mu(n)^{2}}{n} = \frac{6}{\pi^{2}} \log \xi + O(1),
 \end{equation}
we conclude that  
\[ I_{1} = \frac{3}{\pi^{3}} T \log \xi + O\Big(\xi \sqrt{\log T} + T \Big)
\]
for our choice of $c$.  Here we have used the fact that
$$  \sum_{n=1}^{\infty}  \frac{\mu(n)^2}{n^{2c-1}} \leq \zeta(2c-1) \ll \log T.$$

To estimate the contribution from the integral $I_{2}$, we recall the functional equation for the Riemann zeta-function which says that
\begin{equation}\label{eq:fe}
\zeta(s)=\chi(s)\zeta(1\!-\!s) 
\end{equation}
where
\begin{equation*}
\chi(s) = 2^{s}\pi^{s-1}\Gamma(1\!-\!s)\sin\big(\frac{\pi s}{2}\big).
\end{equation*}
Stirling's asymptotic formula for the Gamma-function can be used to show that
\begin{equation}
 \label{eq:negstir}
\big| \chi(\sigma\!+\!it) \big| = \Big(\frac{|t|}{2\pi}\Big)^{1/2-\sigma}\Big(1+O\big(|t|^{-1}\big)\Big)
\end{equation}
uniformly for $-1\leq \sigma\leq 2$ and $|t|\geq 1$.  Combining this estimate and (\ref{neg3}), it follows that, for $T\in\mathcal{T}$, 
$$ \big|\zeta(\sigma+iT)\big|^{-1} \ll T^{\min(\sigma-1/2), 0)}\exp\Big(\frac{A\log T}{\log \log T}\Big)$$
uniformly for $-1\leq \sigma\leq 2$.   In addition, we have 
the trivial bound 
\begin{equation}
  \label{eq:Mxibd}
|M_{\xi}(\sigma+it)| \ll \xi^{1-\sigma}.
\end{equation}
Thus, estimating the integral $I_2$ trivially, we find that 
\[
  I_{2} \ll \exp\Big(\frac{A\log T}{\log \log T}\Big)
  \int_{1-c}^{c} T^{\min(\sigma-1/2), 0)} \xi^{\sigma} d \sigma
  \ll \xi \exp\Big(\frac{A\log T}{\log \log T}\Big).
\]

To bound the contribution from the integral $I_{3}$, we notice that the functional equation for $\zeta(s)$ combined with the estimate in \eqref{eq:negstir} implies that, for $1\leq |t|\leq T$,
$$ \big|\zeta(1\!-\!c\!+\!it)\big|^{-1} \ll |t|^{1/2-c}\big|\zeta(c\!-\!it)\big|^{-1} \ll  |t|^{1/2-c} \zeta(c) \ll  |t|^{-1/2} \log T.$$
It therefore follows that
\begin{equation*}
 I_{3} \ll  \log T \Big(\sum_{n\leq \xi} \frac{|\mu(n)|}{n^{c}} \Big) \int_{1}^{T} t^{-1/2} \  dt
   \ll \sqrt{T} (\log T) \log \xi. 
 \end{equation*}

Finally, since $1/\zeta(s)$ and $\mathcal{M}_{\xi}(1\!-\!s)$ are bounded on the interval $[1-c+i, c+i]$, we find that $I_{4} \ll 1$.  Hence, our combined estimates for $I_{1}, I_{2}, I_{3},$ and $I_{4}$ imply that
$$ M_{1} = \frac{3}{\pi^{3}} T \log \xi+ O\Big( \xi \exp\Big(\frac{A \log T}{\log\log T}\Big) + T\Big).$$
From this and \eqref{eq:Tth}, the estimate in \eqref{eq:m1} follows.


\section{The estimation of $M_{2}$}

We now turn our attention to estimating the sum $M_{2}$.  As before, let $T \geq 4$ and $c =1 + (\log T)^{-1}$.  
Assuming the Riemann Hypothesis, we notice that
\[
   M_{2}= \sum_{0<\gamma\leq T} \mathcal{M}_{\xi}(\rho)\mathcal{M}_{\xi}(1\!-\!\rho). 
\]
Therefore, by the residue theorem, we see that
\begin{align*}
M_{2}& = \frac{1}{2\pi i} \left(\int_{c+i}^{c+iT} + \int_{c+iT}^{1-c+iT}+\int_{1-c+iT}^{1-c+1}+\int_{1-c+i}^{c+i}\right) \  M_{\xi}(s) M_{\xi}(1\!-\!s)  \frac{\zeta'}{\zeta}(s)  \ \! ds  \\
   &   =  J_1 + J_2 + J_3 + J_4,
\end{align*}
say. In order to evaluate the integrals over the horizontal part of the contour we shall impose some extra conditions on $T$. 
Without loss of generality, we may assume that  $T$ satisfies
\begin{equation}
\begin{split}
  \label{eq:zpbd}
  |\gamma-T| &\gg \frac{1}{\log T}   \  \text{ for all ordinates } \gamma \text{ and} \\
  \frac{\zeta'}{\zeta}(\sigma\!+\!iT) &\ll (\log T)^2  \ \text{ uniformly for all } 1-c \le \sigma \le c.
\end{split}
\end{equation}
In each interval of length one such a $T$ exists.  This well-known argument may be found in \cite{Da}, page 108.
Applying \eqref{eq:Mxibd} we find
\[
   \sum_{T < \gamma < T+1} |M_{\xi}(\rho) M_{\xi}(1\!-\!\rho)| \ll \xi (\log T). 
\]
Therefore our choice of $T$ determines $M_2$ up to an error term $O(\xi \log T)$. 
First we estimate the horizontal portions of the contour.  By \eqref{eq:Mxibd} and \eqref{eq:zpbd}, we have
\begin{align*}
  J_2=  \frac{1}{2 \pi} \int_{c}^{1-c}  M_{\xi}(\sigma\!+\!it)M_{\xi}(1\!-\!\sigma\!-\!it) \frac{\zeta'}{\zeta}(\sigma\!+\!it) \ \! d\sigma 
  \ll \xi (\log T)^2. 
\end{align*}
Similarly, it may be shown that $J_4 \ll \xi$.  Next we relate $J_3$ to $J_1$. 
We have
\begin{align*}
  J_3 & = \frac{1}{2 \pi} \int_{T}^{1} M_{\xi}(1\!-\!c\!+\!it)M_{\xi}(c\!-\!it)  \frac{\zeta'}{\zeta}(1\!-\!c\!+\!it) \ \! dt \\
  & = -\frac{1}{2 \pi}\overline{ \int_{1}^{T} M_{\xi}(1\!-\!c\!-\!it)M_{\xi}(c\!+\!it)  \frac{\zeta'}{\zeta}(1\!-\!c\!-\!it) \ \! dt} \
\end{align*}
By differentiating  \eqref{eq:fe}, the functional equation, we find that
\begin{align*}
   -\frac{\zeta'}{\zeta}(1\!-\!c\!-\!it) =  - \frac{\chi'}{\chi}(1\!-\!c\!-\!it) + \frac{\zeta'}{\zeta}(c\!+\!it)
\end{align*}
and hence that
\begin{align*}
  J_3 & = -\frac{1}{2 \pi}\overline{ \int_{1}^{T} M_{\xi}(1\!-\!c\!-\!it)M_{\xi}(c\!+\!it)  \frac{\chi'}{\chi}(1\!-\!c\!-\!it) \ \! dt  } \\
        &\quad  + \frac{1}{2 \pi}\overline{ \int_{1}^{T} M_{\xi}(1\!-\!c\!-\!it)M_{\xi}(c\!+\!it)  \frac{\zeta'}{\zeta}(c\!+\!it) \ \! dt  }.
\end{align*}
By \eqref{eq:fe} and Stirling's formula it can be shown that
\[
    -\frac{\chi'}{\chi}(1\!-\!c\!-\!it) = \log \Big(\frac{|t|}{2 \pi} \Big) (1+ O(|t|^{-1}))
\]
uniformly for $1\leq |t| \leq T$.
By \eqref{eq:Mxibd}, the term $O(|t|^{-1})$ contributes to $J_3$ an amount which is $O(\xi \log T)$ and, hence,
it follows that  
\[
  J_3 =  K + \overline{J_1}+O\big(\xi (\log T)\big)
\]
where
\[
   K = \int_{1}^{T}   \log \Big( \frac{t}{2 \pi} \Big) M_{\xi}(c\!+\!it) M_{\xi}(1\!-\!c\!-\!it) \ \! dt.
\]
Collecting estimates, we deduce that
\begin{equation}
   \label{eq:M2}
   M_2 = K + 2 \Re J_1 + O\big(\xi (\log T)^2\big). 
\end{equation}
To complete our estimation of $M_2$, it remains to evaluate $K$ and then $J_1$. 
Integrating by parts, it follows that 
\begin{align*}
   K & = \frac{1}{2 \pi}
    \log  \Big( \frac{T}{2 \pi} \Big) \int_{1}^{T} M_{\xi}(c\!+\!it) M_{\xi}(1\!-\!c\!-\!it) \ \! dt  \\
    &\quad  - \frac{1}{2 \pi} \int_{1}^{T}   \Big(  \int_{1}^{t} M_{\xi}(c\!+\!iu) M_{\xi}(1\!-\!c\!-\!iu)du  \Big) \frac{dt}{t}. 
\end{align*}
By \eqref{eq:mv}, we have 
\begin{align*}
    \int_{1}^{t} M_{\xi}(c\!+\!iu) M_{\xi}(1\!-\!c\!-\!iu) \ \! du & = (t\!-\!1) \sum_{n \le \xi} \frac{\mu(n)^2}{n}
    + O(\xi \sqrt{\log T})   \\
    & = \frac{6}{\pi^2} t \log \xi +  O(\xi \sqrt{\log T} +t)
\end{align*}
for $t > 1$. Substituting this estimate into the above expression for $K$, we see that
\begin{equation}
\begin{split}
  \label{eq:K}
  K & =   \frac{3}{\pi^3} T  \log  \Big(  \frac{T}{2 \pi} \Big) \log \xi + O(T \log T)  +O(T \log \xi)  \\
    & = \frac{3}{\pi^3} T  \log  \Big( \frac{T}{2 \pi} \Big)  \log \xi + O(T \log T). 
\end{split}
\end{equation}
We finish by evaluating the integral $J_1$ which is similar to the evaluation of the integral $I_1$ in the previous section.
By another application of \eqref{eq:mv}, we find that
\begin{align*}
  J_1 =  - \frac{1}{2 \pi} \int_{1}^{T} \sum_{n=1}^{\infty} \frac{\alpha_n}{n^{c+it}}  \sum_{n \le \xi} \frac{\mu(n)}{n^{1-c-it} } \ \! dt
   & = - \frac{(T\!-\!1)}{2 \pi} \sum_{n \le x} \frac{\alpha_n \mu(n)}{n}  \\
   &\quad + 
   O \Bigg(
  \Big(  \sum_{n=1}^{\infty}  \frac{\alpha_{n}^2}{n^{2c-1}} \Big)^{\frac{1}{2}}
    \Big(  \sum_{n \le \xi} \frac{\mu(n)^2}{n^{1-2c}} \Big)^{\frac{1}{2}}
   \Bigg)
\end{align*}
where the coefficients $\alpha_n$ are defined by
 $$\alpha_n = \sum_{\substack{k\ell=n \\ \ell \le \xi}} \Lambda(k) \mu(\ell).$$
Observe that trivially $|\alpha_n| \le \sum_{u \mid n} \Lambda(u) \le  \log n$.  It follows that the error term in the above expression for $J_1$ is 
 $\ll \zeta''(2c-1)^{\frac{1}{2}} \xi \ll \xi (\log T)^{\frac{3}{2}}$.   Finally, we note that
 \begin{align*}
    \sum_{n \le x} \frac{\alpha_n \mu(n)}{n} 
    & =\sum_{\ell \le x} \frac{\mu(\ell)}{\ell} \sum_{k \le \frac{x}{\ell}} \frac{\Lambda(k) \mu(k\ell)}{k} = \sum_{\ell \leq \xi} \frac{\mu(\ell)}{\ell} \!\!\!\!
 \sum_{{\begin{substack}{ p^{j} \le \xi/\ell
         \\  p \text{ prime}, \ j \ge 0
         }\end{substack}}} \!\!\!\!\!\!
  \frac{\mu(p^{j}\ell) \log p }{p^{j}} \\
   & =  \sum_{\ell \leq \xi} \frac{\mu(\ell)}{\ell}\sum_{p\leq \xi/\ell} \frac{\mu(p\ell) \log p }{p} + O(  \log \xi) \\
   &= - \sum_{\ell \leq \xi} \frac{\mu(\ell)^{2}}{\ell}\sum_{p \le\xi/\ell} \frac{\log p}{p} + O\Big(\log \xi  + \sum_{\ell \le \xi} \frac{1}{\ell} \sum_{p|\ell} \frac{\log p}{p}\Big)
 \end{align*}
since $\mu(p\ell)=-\mu(\ell)$ if $(p,\ell)=1$ and $\mu(p\ell) = 0 = O(1)$ if $p|\ell$. 
The sum in the error term is 
\[
 \sum_{\ell \le \xi} \frac{1}{\ell} \sum_{p|\ell} \frac{\log p}{p} = \sum_{ p \le x} \frac{(\log p)}{p^2} \sum_{\ell' \le \frac{\xi}{p}} \frac{1}{\ell'}
\ll \log \xi. 
\]
Hence, by the elementary result $ \sum_{p \leq \xi } \frac{\log p}{p} = \log \xi + O(1)$, \eqref{eq:neg4}, and partial summation, we deduce that
\begin{equation*}\label{sig1}
\begin{split}
\sum_{n \le x} \frac{\alpha_n \mu(n)}{n}    = - \sum_{l\leq \xi} \frac{\mu(l)^{2}\log(\frac{\xi}{l})}{l} + O(\log \xi ) 
= - \frac{3}{\pi^{2}} (\log \xi)^2  + O(\log \xi).
\end{split}
\end{equation*}
Therefore, combining formulae, we have 
\begin{equation}
   \label{eq:J1}
   J_1  =  - \frac{3}{2 \pi^3} T (\log \xi)^2 + O(T \log T).
\end{equation}
Finally \eqref{eq:M2}, \eqref{eq:K}, and \eqref{eq:J1} imply that
$$M_{2} =  \frac{3}{\pi^{3}} T \log T \log \xi +  \frac{3}{\pi^{3}} T (\log \xi)^{2} + O(T\log T)$$
and, thus, by \eqref{eq:Tth}  we deduce \eqref{eq:m2}.

\end{document}